\documentclass[]{siamltex}
\usepackage{amsmath}
\usepackage{amssymb}
\usepackage{graphicx}

\usepackage{tikz}
\usepackage{pgfplots}
\usepackage[Publish,Externaldir={./}]{externaltikz}
\usepackage{calc}

\usepackage{hyperref}

\newcommand{\dis}{\displaystyle}
\newcommand{\R}{\mathbb{R}} 

\newtheorem{remark}{Remark}

\begin{document}

\title{Kinetic layers and coupling conditions for  macroscopic equations on networks I: the wave equation}

\author{R. Borsche\footnotemark[1] 
     \and  A. Klar\footnotemark[1] \footnotemark[2]}
\footnotetext[1]{Technische Universit\"at Kaiserslautern, Department of Mathematics, Erwin-Schr\"odinger-Stra{\ss}e, 67663 Kaiserslautern, Germany 
  (\{borsche, klar\}@mathematik.uni-kl.de)}
\footnotetext[2]{Fraunhofer ITWM, Fraunhoferplatz 1, 67663 Kaiserslautern, Germany} 
 
\date{}


\maketitle

\begin{abstract}
We consider kinetic and associated macroscopic equations on networks. 
The general approach will be explained in this paper for a linear kinetic BGK model and the corresponding limit for small Knudsen number, which is  the wave equation.
Coupling conditions for the macroscopic equations are derived from the kinetic conditions via an asymptotic analysis near the nodes of the network. 
This analysis leads to the consideration of a fixpoint problem involving the coupled solutions of kinetic half-space problems. 
A new approximate method for the solution of kinetic half-space problems is derived and used for the determination of the coupling conditions. 
Numerical comparisons between the solutions of the macroscopic equation with different coupling conditions and the kinetic solution are presented for the case of tripod and more complicated networks.
\end{abstract}




\section{Introduction}

There have been many attempts  to define coupling conditions for macroscopic partial differential equations on networks including, for example,  drift-diffusion equations, scalar  hyperbolic equations, or hyperbolic systems like  the wave equation or Euler type models, see for example \cite{BGKS14,CC17, BNR14,BHK06a,BHK06b,CHS08,ALM10,EK16,VZ09,BCG10,HKP07,CGP05}. 
In \cite{CGP05,HKP07} coupling conditions for scalar hyperbolic equations on networks are discussed and investigated.
\cite{ALM10,EK16,VZ09} treat the wave equation and general nonlinear hyperbolic systems are considered in \cite{BNR14,BHK06a,BHK06b,CHS08,BCG10}. 
We finally note, that, for example, for hyperbolic systems  on networks 
there are still many unsolved problems, like finding suitable coupling conditions without restricting to subsonic situations.

On the other hand, coupling conditions for kinetic equations on networks have been discussed in a much smaller number of publications,
see \cite{FT15,HM09,BKKP16}. 
In \cite{BKKP16} a first attempt to derive a coupling condition for a
macroscopic equation from the underlying kinetic model has been presented for the case of a kinetic equations for chemotaxis.
In the present  paper, we will present a more general and more accurate procedure.
It is motivated by the classical procedure to find kinetic slip boundary conditions for macroscopic equations based on the analysis of the 
kinetic layer \cite{BSS84,BLP79,G92,G08,UTY03}.
In this work, we will derive coupling conditions for macroscopic equations 
on a network from underlying microscopic or kinetic models via an asymptotic analysis of the situation near the nodes. 
To explain the basic approach we concentrate on a simple linear BGK-type kinetic model with the linear wave equation as the associated macroscopic model. 
More complicated problems and, in particular, nonlinear models will be discussed in future work. 

The paper is organized in the following way. 
In section \ref{equations} we discuss the kinetic and macroscopic equations and the boundary and coupling conditions for these equations.
In section \ref{layeranalysis}  kinetic boundary layers are discussed, as well as   an 
asymptotic analysis of the kinetic equations near the nodes. 
This leads to an abstract formulation of the coupling conditions for the macroscopic equations at the nodes based on a fix-point problem involving kinetic half-space equations. 
In the following section \ref{halfspacemarshak} an approximate coupling condition is derived based on the so called Maxwell approximation of the half-space problem. 
A  refined method to determine the solution of the half space problems is derived, compared to previous approximate solution methods for half-space problems and applied to the problem of finding accurate coupling conditions for the macroscopic equations in sections \ref{half moment coupling}. 
Moreover, the macroscopic equations on the network with the different coupling conditions are numerically compared to each other and to the full solutions of the kinetic equations on the network  in section \ref{Numerical results}.
The results show the very good approximation of the underlying kinetic model by the macroscopic model with the new coupling conditions.

\section{Equations and boundary and coupling conditions}
\label{equations}

We consider  the  linear kinetic  BGK model in 1D for  
  $x \in \mathbb{R}, v  \in \R$.
\begin{align}\label{bgk}
\begin{aligned}
\partial_t f + v \partial_x f = -\frac{1}{\epsilon} \left(f-\left(\rho + \frac{v}{ a^2}  q \right) M (v)\right) 
\end{aligned}
\end{align}
with the Maxwellian
$$
M (v) = 
\frac{1}{\sqrt{2 \pi a^2}} \exp \left(\frac{ - \vert v\vert^2}{2 a^2} \right),
$$
where
$$
\rho_\epsilon = \int_{-\infty}^\infty f(v) dv ,\qquad  q_\epsilon = \int_{-\infty}^\infty v f(v) dv .
$$
The associated macroscopic equation  for $\epsilon \rightarrow 0$ is the wave equation
\begin{align}\label{euler}
\begin{aligned}
\partial_t \rho + \partial_x q =0\\
\partial_t q + \partial_x (a^2 \rho ) =0\ .
\end{aligned}
\end{align}

To illustrate the influence of the underlying kinetic model on the coupling conditions for the macroscopic equations, we also consider the following equation with  bounded velocity space  $ v  \in [-1,1]$ 
\begin{align}\label{bv}
\begin{aligned}
\partial_t f + v \partial_x f = -\frac{1}{\epsilon} \left(  f-  (\rho_\epsilon+ \frac{v}{a^2}  q_\epsilon ) \frac{1}{2} \right)
\end{aligned}
\end{align}
with $a^2 = \frac{1}{3}$ or 
\begin{align}\label{bgkbounded}
\begin{aligned}
\partial_t f + v \partial_x f = -\frac{1}{\epsilon} \left(  f-  (\frac{\rho_\epsilon}{2}+ \frac{3}{2} v q_\epsilon ) \right)
\end{aligned}
\end{align}
with
$$
\rho_\epsilon = \int_{-1}^1 f(v) dv ,\quad  q_\epsilon = \int_{-1}^1 v f(v) dv \ .
$$
The associated macroscopic equation is again \eqref{euler} 
with
$a^2 = \frac{1}{3}$.
\subsection{Boundary conditions}
For $x \in [0,b]$ we prescribe for the kinetic  equation  
$$
 f(0,v), v > 0, \; f(b,v), v<0.
$$
For \eqref{euler} the boundary conditions are given in characteristic variables \cite{T09}.
The corresponding Riemann Invariants are
\begin{align}\label{eq:RiemannInvariants}
r_{1/2} = q   \mp a \rho.
\end{align}
As boundary data the value of $ q + a\rho $
at the left boundary and $q - a\rho$ at the right boundary are prescribed.

\subsection{Coupling conditions}
\begin{figure}[h!]
	\begin{center}
		\begin{tikzpicture}
		\draw[->] (0,0)--(-1.4,0) node[left]{$1$};
		\draw[->] (0,0)--(1,-1) node[right]{$2$};
		\draw[->] (0,0)--(1,1) node[right]{$3$};
		\node[fill=black,circle] at (0,0){};
		\end{tikzpicture} 
	\end{center}
	\caption{Node connecting three edges and orientation of the edges.}
	\label{fig:tripod}
\end{figure}
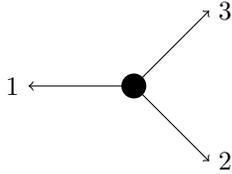
If these equations are considered on a network, it is sufficient to study a single coupling point.
At each node so called coupling conditions are required. 
In the following we consider a node connecting $n$ edges, which are oriented away from the node, as in Figure \ref{fig:tripod}.
Each edge $i$ is parametrized by the interval $[0,b_i]$ and the kinetic and macroscopic quantities are denoted by
$f^i$ and $\rho^i, q^i$ respectively.
A possible choice of coupling conditions for the kinetic problem are given by 
\begin{align}
\label{kincoup}
f^i(0,v) = \sum_{j=1}^n c_{ij}f^j (0,-v), v >0,
\end{align}
compare \cite{BKKP16}.
The total mass in the system is conserved, if 
\begin{align}\label{eq:conservative_coupling_f}
\sum_{i=1}^n c_{ij} = 1
\end{align}
 holds.
In the following we use the vector notation 
\begin{align*}
	f^+ = C f^-, v >0\ ,
\end{align*}
where $f^+=(f^1(0,v),\dots,f^n(0,v))$ and $f^-=(f^1(0,-v),\dots,f^n(0,-v))$.
The coupling conditions for the macroscopic quantities are
conditions on the characteristic variables
\begin{align*}
r_2^i (0) = q^i(0) + a \rho^i(0)
\end{align*}
using the given values of 
\begin{align*}
r_1^i (0) = q^i(0) - a \rho^i(0)\ .
\end{align*}
We refer to \cite{CHS08,G10} for systems of macroscopic equations on networks.
In the following we will derive, via asymptotic analysis, 
macroscopic coupling conditions from the kinetic coupling conditions (\ref{kincoup}).

\section{Boundary and coupling conditions for macroscopic equations via kinetic layer analysis}
\label{layeranalysis}
\subsection{Boundary conditions}
A kinetic layer analysis, see \cite{BSS84,C69,CGS88,G10}, at the left boundary of the interval $[0,b]$, i.e.  a rescaling of the spatial variable in equation (\ref{bgk}) with $\epsilon$, gives   to first order in $\epsilon$  the  following stationary kinetic  half space  problem for 
$x \in [0,\infty]$
\begin{align}\label{bgkhalfspace}
\begin{aligned}
 v \partial_x \varphi = - \left( \varphi -\left(\rho + \frac{v}{ a^2} q \right) M (v) \right)\ ,
\end{aligned}
\end{align}
where $\rho$ and $q$ are the zeroth and first moments of  $\varphi$. 
At $x=0$ the boundary conditions for the  half space problem are
$$
\varphi(0,v) =  k(v) = f(0,v), v > 0\ .
$$

On the left side, i.e. at $x=\infty$,  a condition is prescribing 
an arbitrary linear combination of the invariants of the half-space problem
$< v \varphi > $
and $< v^2\varphi > $ .
Here, and in the following we use the notation $< \varphi  > = \int \varphi dv $ and  $< \varphi  >_+ = \int_{v>0} \varphi dv $
or $< \varphi  >_- = \int_{v<0} \varphi dv $.
We use the values of the first Riemann Invariant \eqref{eq:RiemannInvariants} $r_1 = q -a  \rho$ of the macroscopic system \eqref{euler} to fix 
$$< \left(v - \frac{v^2}{ a } \right)\varphi> = r_1.
$$

The boundary condition for \eqref{euler} is obtained by 
determining $r_2$ from  the asymptotic solution
of the half space problem and setting $r_2 = q_\infty +a  \rho_\infty$.
The values $ \rho_\infty$ and $ q_\infty  $ are the macroscopic quantities associated to the solution of the half-space problem at infinity, which has the form
$$
\varphi(\infty,v) = \left(\rho_\infty+ \frac{v}{a^2} q_\infty \right)M(v)\ .
$$

\begin{figure}[h!]
	\begin{center}
		\begin{tikzpicture}
		\draw[->] (-1,0)--(1.5,0) node[right]{$f$};
		\draw[->] (0,-2)--(0,2) node[left]{$v$};
		\draw (-0.5,0) node[below]{$x=0$};
		\draw[->] (-0.2,1)--(0.2,1);
		\draw (-0.2,1.5) node[left]{$k(v)$};
		\draw[->] (0.2,-1)--(-0.2,-1);
		\draw (-0.2,-1.5) node[left]{$A[k](v)$};
		
		\draw[domain=0:2,gray] plot ({exp(-((\x)^2)/1)},\x); 
		\draw[domain=0:2,gray,dashed] plot ({exp(-((\x)^2)/1)},-\x); 
		
		\end{tikzpicture} 
		\qquad 
		\begin{tikzpicture}
		\draw[->] (-1,0)--(5.5,0) node[right]{$x$};
		\draw[->] (0,-0.25)--(0,3.2) node[left]{$\rho$};
		\draw[-] (4.25,0.25)--(4.25,-0.25) node[below]{$x_\infty$};
		\draw (5.,1) node[below]{$(\rho,q)$};
		\draw (3.5,1) node[below]{$(\rho_\infty,q_\infty)$};
		\draw (0,-0.25) node[below]{$x=0$};
		
		\draw[domain=0:5,gray] plot (\x,{2*exp(-((\x)^2)/0.5)+1}); 
		
		\end{tikzpicture}
	\end{center}
	\caption{Illustration of the half space problem}
	\label{fig:HalfSpace}
\end{figure}
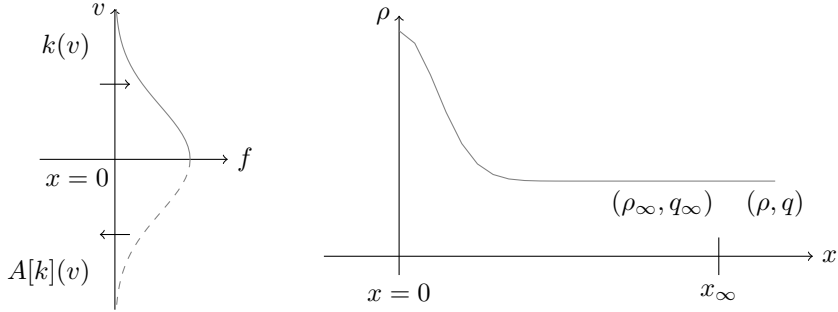

The solution of the half space problem is also used to determine the outgoing distribution
$$
A[k](v) = f(0,v)  = \varphi(0,v), v<0\ ,
$$
where the notation $A$ is used for the so called Albedo operator of the  half space problem. 
The structure of a half space problem is illustrated in Figure \ref{fig:HalfSpace}.

\subsection{Coupling conditions}

We use the corresponding procedure to determine the coupling conditions for the macroscopic equations.
\begin{figure}[h!]
	\begin{center}
		\begin{tikzpicture}
		\draw[->] (0,0)--(-2*1.4,0) node[left]{$1$};
		\draw[->] (0,0)--(2*1,-1*2) node[right]{$2$};
		\draw[->] (0,0)--(2*1,1*2) node[right]{$3$};
		\draw (-2*1.1,0.1)--(-2*1.1,-0.1) node[below]{$x_\infty$};
		\draw (1.6,1.8)--(1.8,1.6) node[right]{$x_\infty$};
		\draw (1.8,-1.6)--(1.6,-1.8) node[left]{$x_\infty$};
		\draw (-0.2,0.2) node[above]{$C$};
		\end{tikzpicture} 
	\end{center}
	\caption{Illustration of a three way junction. Three connected half spaces.}
		\label{fig:junction}
\end{figure}
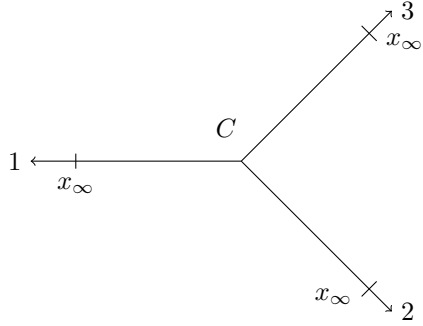

Starting from the kinetic coupling conditions
\begin{align*}
f^i(0,v) = \sum_{j=1}^n c_{ij}f^j (0,-v), v >0
\end{align*}
we determine
the coupling conditions for the macroscopic equations in the following way.
We use the kinetic  coupling conditions to obtain conditions on the in- and outgoing solutions of the half space problems on the different arcs.
That means 
\begin{align*}
\varphi^i(0,v)
= \sum_{j=1}^n c_{ij} \varphi^j(0,-v) , v>0
\end{align*}
or, if we denote the ingoing function of the half-space problem on arc $i$ by
$k^i (v),v>0$  and the outgoing solution by $A^i [k^i] (v), v<0$
\begin{align}
\label{fixpoint}
 k^i(v)
 = \sum_{j=1}^n c_{ij}  A^j [ k^j] (-v) , v  >0\ .
\end{align}

This is a fix point equation for $k^i, i= 1 \ldots,n$.
Additional conditions are needed to solve the half-space problems, i.e.

$$< \left( v- \frac{v^2}{a } \right)\varphi^i>= r^i_1 
$$
with $r^i_{1} =  q^i -    a \rho^i$.
The coupling conditions for the wave equation  are  conditions on the outgoing characteristic variables at $x=0$. We define 
\begin{align*}
r^i_{2} (0) = 
 q^i_\infty[k^i]+ a \rho^i_\infty[k^i].
\end{align*}

The main task is now to find tractable expressions for these 
coupling conditions. In the next section we discuss the 
so called Maxwell method to solve the half-space problem approximately and the resulting approximate coupling conditions.
In further sections a new refined method to determine 
the solution of the half space problems based on half moment equations is derived and applied to the problem of finding accurate coupling conditions. 
Finally, the results are compared to full solutions of the kinetic equations on the network and to other approximate solution methods for the half-space problem.

\section{Approximate solution  of the half space problem via half-fluxes and approximate coupling conditions}
\label{halfspacemarshak}

To solve kinetic  half space problem approximately a variety of different methods can be found in the literature. 
Approaches via a Galerkin method can be found in \cite{Coron,LLS16,LLS162}. 
Approximate  methods to determine only the  asymptotic states and outgoing distributions can be found in \cite{GK95,LF67,LF81}.
For the determination of the macroscopic  boundary conditions we use in this section  a simple approximation, the so called Maxwell approximation, see \cite{M,M47}. 

\subsection{Approximate solution of the half space problem}
\label{boundmaxwell}

We use  the equality of half-fluxes
\begin{align*}
<v k(v) >_+  
= <v \left(\rho_\infty + \frac{v}{ a^2} q_\infty \right) M >_+
\end{align*}
to obtain one condition. 
Together with the first Riemann Invariant
\begin{align*}
 q_\infty -a\rho_\infty
 = q -a \rho = r_1(0) = C = \mbox{const.}
\end{align*}
this determines $\rho_\infty $ and $q_\infty$.
This can be rewritten as
\begin{align*}
\left(\begin{array}{cc} a \frac{1}{\sqrt{2 \pi}}  & \frac{1}{2} \\ -a  & 1 \end{array} \right)
\left(\begin{array}{c}  \rho_\infty\\ q_\infty \end{array} \right)
=
\left(\begin{array}{c} <vk>_+\\ C\end{array} \right)\ 
\end{align*}
with the solution
\begin{align*}
\left(\begin{array}{c}  \rho_\infty\\ q_\infty \end{array} \right)
=
\frac{1}{a\left( \frac{1}{\sqrt{2 \pi}}- \frac{1}{2}\right)}\left(\begin{array}{c} <vk>_+- \frac{C}{2}\\  a <vk>_+ + \frac{aC}{\sqrt{2 \pi}}\end{array} \right)
\ .
\end{align*}
The outgoing distribution is  approximated by
$$
\varphi(0,v)  =  \left( \rho_\infty + \frac{v}{ a^2}  q_\infty \right) M (v ) , v <0.
$$

\begin{remark}
For bounded velocity space and equation \eqref{bv} we have
\begin{align}\label{eq:Maxwell_bounded}
\left(\begin{array}{cc} \frac{1}{4}  & \frac{1}{2} \\ -a  & 1 \end{array} \right)
\left(\begin{array}{c}  \rho_\infty\\ q_\infty \end{array} \right)
=
\left(\begin{array}{c} <vk>_+\\ C\end{array} \right)
\ .
\end{align}
\end{remark}

\subsection{Approximate coupling conditions}

The fix-point problem (\ref{fixpoint}) is in the present case 
approximated by the problem
\begin{align*}
k^i(v) = \sum_{j=1}^n c_{ij} \left( \rho_\infty^j [k^j]- \frac{v}{ a^2}  q_\infty^j [k^j]
\right) M (-v) ,  v >0.
\end{align*}
Thus, the asymptotic values are determined by
\begin{align}\label{eq:halfspace_Maxwell}
\left(\begin{array}{cc} a \frac{1}{\sqrt{2 \pi}}  & \frac{1}{2} \\ -a  & 1 \end{array} \right)
 \left(\begin{array}{c} \rho^i_\infty \\ q^i_\infty \end{array} \right)
=
\left(\begin{array}{c} \sum_{j} c_{ij} \left( \rho_\infty^j \frac{a}{\sqrt{2 \pi}} - \frac{ q_\infty^j}{ 2}  \right) \\  r^i_1 (0) \end{array} \right).
\end{align}
These values can be assigned to the states of the wave equation \eqref{euler} at the coupling point by
$$
r^i_{2} (0) = q^i_\infty +a \rho^i_\infty .\\
$$
Summing the first equation of \eqref{eq:halfspace_Maxwell} for all $i=1,\dots,n$ and using the conservation property of the kinetic coupling conditions \eqref{eq:conservative_coupling_f} directly yields the conservation of mass in the macroscopic variables
\begin{align}\label{eq:sum_q_0}
\sum_{i=1}^n q_\infty^i =0.
\end{align}
In the special case of a  uniform node, i.e. $c_{ij} = \frac{1}{n-1}$ for $i \neq j$ and $c_{ij}=0$ for $i =j$  we have
\begin{align*}
 a \frac{1}{\sqrt{2 \pi}} \rho^i_\infty +   \frac{1}{2}  q^i_\infty 
=\frac{1}{n-1}\sum_{j\ne i} 
 \left(\frac{a}{ \sqrt{2 \pi}}  \rho^j_\infty  -   \frac{q^j_\infty }{ 2} \right). 
\end{align*}
Multiplying by $\sqrt{2 \pi}$ and adding $\frac{\sqrt{2 \pi}}{n}$ times \eqref{eq:sum_q_0} we obtain
\begin{align*}
a  \rho^i_\infty +   \frac{\sqrt{2 \pi}}{2} \frac{n-2}{n} q^i_\infty 
=\frac{1}{n-1}\sum_{j\ne i} 
\left(a  \rho^j_\infty  +\frac{\sqrt{2 \pi}}{2} \frac{n-2}{n}   q^j_\infty \right)\qquad i = 1,\dots,n. 
\end{align*}
These $n$ equations are not linearly independent and they can be reformulated as
\begin{align*}
 a \rho^i_\infty +   \frac{\sqrt{2 \pi}}{2} \frac{n-2}{n} q^i_\infty 
=
a   \rho^j_\infty  +   \frac{\sqrt{2 \pi}}{2}\frac{n-2}{n} q^j_\infty  
 \qquad i,j=1,\dots,n\ , i\neq j. 
\end{align*}
Thus we have found a macroscopic invariant at the coupling point, i.e.
\begin{align}
\label{maxwellunbounded}
 \rho +   \frac{\sqrt{2 \pi}}{2 a } \frac{n-2}{n} q. 
\end{align}

\begin{remark}
For bounded velocity space we have
\begin{align*}
  \frac{1}{4} \rho^i_\infty +   \frac{1}{2}  q^i_\infty 
=\sum_{j} c_{ij} \left(
 \frac{1}{ 4}  \rho^j_\infty -   \frac{q^j_\infty }{ 2}  \right).
\end{align*}
For a uniform node this leads to
\begin{align*}
  \rho^i_\infty +   \frac{2(n-2)}{n}  q^i_\infty 
=
 \rho^j_\infty +   \frac{2(n-2)}{n} q^j_\infty  
\end{align*}
and to the invariant
\begin{align}\label{eq:Invariants_Maxwell}
\rho_\infty +   \frac{2(n-2)}{n}  q_\infty .
\end{align}
For example for $n=3$, we obtain for  the above invariant \eqref{maxwellunbounded}  a factor 
$\frac{\sqrt{6 \pi}}{6} \sim 0.72360$ for $a^2 = \frac{1}{3}$ in contrast to the value $\frac{2}{3}$ for the bounded velocity case. 
Thus although the kinetic coupling conditions and the macroscopic equations \eqref{euler} are identical we obtain different coupling conditions for a model with bounded velocity, since the half space problems are different.
\end{remark}

\begin{remark}
 Another direct approach to obtain coupling conditions for the wave equation has been used in \cite{BKKP16}.
 Using a full moment approximation of the distribution function in the case of bounded velocities, i.e. using
 $$
 f^\pm (v) = \frac{\rho}{2} \pm \frac{3}{2} v q,\quad  v\in [0,1]
 $$
 in the kinetic coupling conditions  (\ref{kincoup}), one obtains after integration over the positive velocities
 \begin{align*}
  \rho^i_\infty  + \frac{3}{2} q^i_\infty =    \sum_j c_{ij}\left(\rho^{j}_\infty - \frac{3}{2} q^{j}_\infty \right)  \ .
  \end{align*} 
 This leads, for a uniform node, to the invariant
   \begin{align}\label{eq:invariant_fullmoment}
      \rho_\infty   + \frac{3(n-2)}{2n} q_\infty \ .
      \end{align} 
      In the case of unbounded velocities we obtain 
      \begin{align*}
        \rho^i_\infty  + \frac{\sqrt{2}}{ a \sqrt{\pi}} q^i_\infty =    \sum_j c_{ij}\left(\rho^{j}_\infty - \frac{\sqrt{2}}{ a \sqrt{\pi}} q^{j}_\infty \right)  \ .
        \end{align*} 
      and the invariant
         \begin{align}\label{eq:invariant_fullmomentunbounded}
            \rho_\infty   + \frac{\sqrt{2}}{ a \sqrt{\pi}} \frac{(n-2)}{n} q_\infty \ .
            \end{align}
  \end{remark}

\section{Coupling conditions via  approximation by half-moment equations}
\label{half moment coupling}
In this section we develop a refined method to solve the half-space problem based on a half moment approximation and use it for the derivation of refined coupling conditions.
We begin with the case of equation \eqref{bv}.
For further approaches to the approximate solution of half space problems we refer to \cite{GK95,LF67,ST81}.

\subsection{The case of a bounded velocity domain}
\label{boundedvel}
Consider  the linear BGK model with bounded velocities \eqref{bgkbounded} and the corresponding limit equation \eqref{euler}.

\subsubsection{Half moment equations} 
\label{boundhalfmoment}

We determine a half moment approximation for the half-space solution, compare \cite{BKKP16}.
We define
\begin{equation}
\label{eq:def_halfmoments}
\begin{aligned}
\rho_\epsilon^- &=\dis\int_{-1}^0 f(v)dv\ , \qquad
& \rho_\epsilon^+ &=\dis\int_{0}^1f(v)dv\ ,
\\
 q_\epsilon^- &=\dis\int_{-1}^0 vf(v)dv \ ,
& q_\epsilon^+ &=\dis\int_{0}^1vf(v)dv \ .
\end{aligned}
\end{equation}
As  closure assumption we use the following approximation of the distribution function $f$ by  affine linear functions in $v$ to determine half-moment equations, see \cite{BKKP16} and references therein.
\begin{align*}
f(v) = a^+ + vb^+,\ v\geq 0&&\text{and}&&
f(v) = a^- + vb^-,\ v\leq 0
\ .
\end{align*}
One obtains
\begin{align*}
\begin{aligned}
\rho_\epsilon^- &= a^- -\dis\frac{1}{2}b^-
\ ,\quad&
\rho_\epsilon^+ &=a^+ + \dis\frac{1}{2}b^+\ , \\[0.3cm]
q_\epsilon^- &= -\dis\frac{1}{2}a^- + \dis\frac{1}{3}b^- 
\ ,\quad &
q_\epsilon^+ &= \dis\frac{1}{2}a^+ + \dis\frac{1}{3}b^+\ 
\end{aligned}
\end{align*}
and 
\begin{align*}
\dis\int_0^1 v^2f(v)dv = -\dis\frac{1}{6}\rho_\epsilon^+ + q_\epsilon^+ \ ,
&&
\dis\int_{-1}^0 v^2f(v)dv = -\dis\frac{1}{6}\rho_\epsilon^- - q_\epsilon^-\ .
\end{align*}
Finally, integrating the kinetic equation with respect to the corresponding half-spaces, we get the half-moment approximation of the  kinetic equation  as
\begin{equation}\label{equ:sc half-moment system}
\left\{
\begin{array}{lcl}
\partial_t\rho_\epsilon^+ + \partial_xq_\epsilon^+ &=& -\dis\frac{1}{\epsilon}\left(\rho_\epsilon^+ - \left(\dis\frac{\rho_\epsilon }{2} + \frac{3 q_\epsilon}{4}\right) \right)\\[0.3cm]
\partial_tq_\epsilon^+ + \partial_x\left(-\dis\frac{1}{6}\rho_\epsilon^+ + q_\epsilon^+\right) &=& -\dis\frac{1}{\epsilon}\left(q_\epsilon^+ -\left(\dis\frac{\rho_\epsilon }{4} + \frac{ q_\epsilon}{2}\right)\right) \\[0.3cm]
\partial_t\rho_\epsilon^- + \partial_xq_\epsilon^- &=& -\dis\frac{1}{\epsilon}\left(\rho_\epsilon^- -
\left(\dis\frac{\rho_\epsilon }{2} - \frac{3 q_\epsilon}{4}\right)
\right)\\[0.3cm]
\partial_tq_\epsilon^- + \partial_x\left( -\dis\frac{1}{6}\rho_\epsilon^- - q_\epsilon^- \right) &=& -\dis\frac{1}{\epsilon}\left(q_\epsilon^- - \left(- \dis\frac{\rho_\epsilon }{4} + \frac{ q_\epsilon}{2}\right)\right) .\\
\end{array}
\right.
\end{equation}
Introducing the even-odd variables
\begin{align*}
\begin{aligned}
\rho_\epsilon &= \rho_\epsilon^+ + \rho_\epsilon^-\ ,
& \qquad
\hat{\rho}_\epsilon &= \rho_\epsilon^+ - \rho_\epsilon^-\ ,
\\
q_\epsilon &= q_\epsilon^+ + q_\epsilon^-\ , 
&
\hat{q}_\epsilon &= q_\epsilon^+ -q_\epsilon^-\ ,
\end{aligned}
\end{align*}
we can rewrite the system as 
\begin{equation*}
\left\{
\begin{array}{lcl}
\partial_t\rho_\epsilon + \partial_xq_\epsilon &=& 0\\[0.3cm]
\partial_tq_\epsilon + \partial_x\left(-\dis\frac{1}{6}\rho_\epsilon + \hat{q}_\epsilon\right) &=&0\\[0.3cm]
\partial_t\hat{\rho}_\epsilon + \partial_x\hat{q}_\epsilon &=&  -\dis\frac{1}{\epsilon}\left(\hat{\rho}_\epsilon -\dis\frac{3}{2}q_\epsilon\right)\\[0.3cm]
\partial_t\hat{q}_\epsilon + \partial_x\left(-\dis\frac{1}{6}\hat{\rho}_\epsilon + q_\epsilon\right) &=& -\dis\frac{1}{\epsilon}\left(\hat{q}_\epsilon-\dis\frac{\rho_\epsilon}{2}\right)\ .
\end{array}
\right.
\end{equation*}
Obviously, the  half-moment model  has again the wave equation \eqref{euler} as  macroscopic limit as $\epsilon$ goes to $0$.

\subsubsection{Half space-half moment problem}

Rescaling the spatial variable in the half-moment problem with $\epsilon$ one obtains the following half-space problem for 
$x \in \R^+$

\begin{equation*}
\left\{
\begin{array}{lcl}
 \partial_xq^+ &=& -\left(\rho^+ - \left(\dis\frac{\rho }{2} + \frac{3 q}{4}\right) \right)\\[0.3cm]
 \partial_x\left(-\dis\frac{1}{6}\rho^+ + q^+\right) &=& -\left(q^+ -\left(\dis\frac{\rho }{4} + \frac{ q}{2}\right)\right) \\[0.3cm]
 \partial_xq^- &=& -\left(\rho^- -
\left(\dis\frac{\rho }{2} - \frac{3 q}{4}\right)
\right)\\[0.3cm]
 \partial_x\left( -\dis\frac{1}{6}\rho^- - q^- \right) &=& -\left(q^- - \left(- \dis\frac{\rho }{4} + \frac{ q}{2}\right)\right) \\
\end{array}
\right.
\end{equation*}
or
\begin{equation}\label{half-moment half-space}
\left\{
\begin{array}{lcl}
 \partial_xq &=& 0\\[0.3cm]
\partial_x\left(-\dis\frac{1}{6}\rho + \hat{q}\right) &=&0\\[0.3cm]
 \partial_x\hat{q} &=&  -\left(\hat{\rho} -\dis\frac{3}{2}q\right)\\[0.3cm]
 \partial_x\left(-\dis\frac{1}{6}\hat{\rho} + q\right) &=& -\left(\hat{q}-\dis\frac{\rho}{2}\right)\ .
\end{array}
\right.
\end{equation}
We  have to provide boundary conditions for $\rho^+(0)$ and $q^+ (0)$, as well as a condition at $x=\infty$
$$
q_\infty - a \rho_\infty = r_1 (0) = C.
$$
Then, the half space problem can be solved explicitly.
We determine a solution up to 3 constants which will be fixed with the above 3 conditions.
First, we observe, that we have 2 invariants
\begin{align}\nonumber
q &= C_1 = const \\
\label{eq:halfmoment_C2}
-\frac{\rho}{6} + \hat q &= C_2
\end{align}
From the last equation in \eqref{half-moment half-space} we can deduce that at $x = \infty$
$$
\hat q_\infty = \frac{\rho_\infty}{2}\ .
$$
Combining this with \eqref{eq:halfmoment_C2} gives 
$ \rho_\infty = 3 C_2$ or $\hat q_\infty = \frac{3 C_2}{2} $.
From the third equation of \eqref{half-moment half-space} we obtain 
$\hat \rho_\infty = \frac{3q}{2} = \frac{3 C_1}{2}$.
This simplifies \eqref{half-moment half-space} to 
\begin{equation*}
\left\{
\begin{array}{rcl}
 q&=&C_1\\
 \rho &=& 6 \hat q - 6 C_2\\
 \partial_x\hat{q} &=&  -\left(\hat{\rho} -\dis\frac{3}{2} C_1\right)\\[0.3cm]
 \partial_x\left(-\dis\frac{1}{6}\hat{\rho}\right) &=& -\left(-2 \hat{q} +  3 C_2  \right)\ .
\end{array}
\right.
\end{equation*}
The ODEs for $\hat \rho$ and  $\hat q$ have the solutions 
\begin{align*}
\hat \rho= \gamma \exp(-\frac{2x}{a})   + \hat \gamma \exp(\frac{2x}{a}) + \frac{3 C_1}{2}\\
\hat q = \frac{a}{2}\gamma \exp(-\frac{2x}{a})   -  \frac{a}{2} \hat \gamma \exp(\frac{2x}{a}) +\frac{3 C_2}{2} 
\end{align*}
Since we are looking only for bounded solutions we are left with
\begin{align*}
\hat \rho= \gamma \exp(-\frac{2x}{a})    + \frac{3 C_1}{2}\\
\hat q = \frac{a}{2}\gamma \exp(-\frac{2x}{a})   +\frac{3 C_2}{2} \\
q=C_1\\
\rho = 3 a \gamma \exp(-\frac{2x}{a})   + 3  C_2 .  
\end{align*}
The three  parameters are  fixed with the 3 conditions mentioned above.
At $x=0$ inflow data is given
\begin{align*}
\frac{1}{2}\left( q(0)+\hat q (0) \right) = q_+ (0)\\
\frac{1}{2}\left( \rho (0)+\hat \rho (0) \right) = \rho_+ (0) 
\end{align*}
and the Riemann Invariant at $x=\infty$ gives
\begin{align}\label{eq:Riemann_INvariant_infty}
	q_\infty -a \rho_\infty = C.
\end{align}
Inserting the above determined solution we obtain
\begin{align*}
\frac{1}{2}\left( C_1 +\frac{a}{2}\gamma   +\frac{3 C_2}{2} \right) &= q_+ (0)\\
\frac{1}{2}\left( 3 a \gamma   + 3  C_2 +\gamma     + \frac{3 C_1}{2} \right) &= \rho_+ (0)\ ,
\end{align*}
 which can be rewritten in terms of the asymptotic values
\begin{align}
\label{eq:halfMoment_rqp}
\begin{aligned}
\frac{q_\infty}{2}   +\frac{\rho_\infty}{4} +\frac{a}{4}\gamma  &= q_+ (0)\\
 \frac{3 q_\infty}{4} +  \frac{\rho_\infty}{2}       + \left(\frac{3 a+1}{2}\right) \gamma  &= \rho_+ (0) \ .
\end{aligned}
\end{align} 
  Together with the condition at infinity \eqref{eq:Riemann_INvariant_infty}, this determines  the   asymptotic values $q_\infty,\rho_\infty$ and $\gamma$.
 The   outgoing quantities $\rho_-(0), q_-(0)$ are then determined by 
 \begin{align*}
 \frac{q_\infty}{2}   - \frac{ \rho_\infty}{4} - \frac{a}{4}\gamma  &= q_-(0)\\
  - \frac{3 q_\infty}{4}  +    \frac{1 }{2}  \rho_\infty     + \left(\frac{3 a-1 }{2}\right) \gamma  &= \rho_- (0) \ .
 \end{align*} 
 
 \begin{remark}
With the Maxwell approximation for the half-moment problem in \eqref{eq:Maxwell_bounded} the asymptotic states are determined by 
 \begin{align*}
 q_+ (0) = q_+(\infty)  = \frac{1}{2} (q_\infty + \hat q_\infty)
 = \frac{q_\infty}{2} + \frac{\rho_\infty}{4}
 \end{align*} 
 and the condition at infinity \eqref{eq:Riemann_INvariant_infty}.
 The outgoing quantities are 
 \begin{align*}
  q_- (0) &= q_-(\infty) = - \frac{\rho_\infty}{4} + \frac{ q_\infty}{2}\\
  \rho_-(0) &= \rho_-(\infty)
  = \frac{\rho_\infty}{2} -  \frac{3 q_\infty}{4}\ .
  \end{align*} 
\end{remark}
\subsubsection{The extrapolation length}
 
 To estimate the accuracy of our method, we consider the classical problem of determining the so-called extrapolation length \cite{LF67,ST81,GK95}.
 
 For $x \in \mathbb{R}^+, v  \in [-1,1]$ we consider the half space 
   equation 
  \begin{align*}
  v \partial_x f = - \left(  f-  (\frac{\rho}{2}+ \frac{3}{2} v q )\right)
  \end{align*}
  with
  $\int_{-1}^1 v f dv = q = 0$ and $f (0,v) = v , v >0$.
  That means, we consider
  \begin{align*}
  v \partial_x f = - \left(  f-  \frac{\rho}{2} \right)\ .
  \end{align*}
  The extrapolation length is the value of $ \lambda_\infty = f(\infty,v) = \frac{\rho_\infty}{2}$.
  The Maxwell approximation  gives $\lambda_\infty =   \frac{2}{3} $.
  The above half moment approximation gives
  \begin{align*}
  \frac{\rho_\infty}{4} +\frac{a}{4}\gamma  = q_+ (0) = \frac{1}{3}\\
   \frac{\rho_\infty}{2}       + (\frac{3 a+1}{2}) \gamma  = \rho_+ (0) = \frac{1}{2} .
  \end{align*} 
  This leads to 
   \begin{align*}
   \rho_\infty    = \frac{  9 a+ 4}{6 a+3 }
   \end{align*} 
   and with $a^2 =\frac{1}{3}$ we obtain
   \begin{align*}
   \rho_\infty    = \frac{  3 \sqrt{3} + 4}{2 \sqrt{3}+3 } .
   \end{align*} 
   Thus, the extrapolation length is approximated as $\lambda_\infty \sim 0.7113$.
   The exact value computed from a spectral method is $0.7104$,
   see \cite{Coron,LLS16}. 
   This yields an error for the above half-moment method of approximately $0.1 \%$.
   In contrast, the Maxwell approximation gives $0.6666$, which is an error of $6.1 \% $. 
   The variational method in \cite{LF67,GK95} gives $0.7083$, which is an error of $0.3\%$.

\subsubsection{Half-moment coupling conditions}

In this subsection we determine the coupling conditions on the basis of the half-moment approximation of the half-space problem.
Multiplying with $v$ and integrating the kinetic coupling conditions \eqref{kincoup} with respect to the  positive and negative half moments  gives
\begin{align}
\label{halfcoupling}
q^i_+ (0) = - \sum_{ j=1}^n  c_{ij} q_-^{j} (0) \ .
\end{align}
Inserting the half moment approximations \eqref{eq:halfMoment_rqp} yields
\begin{align*}
\frac{q^i_\infty}{2}   +\frac{\rho^i_\infty}{4} +\frac{a}{4}\gamma^i  
=  \sum_{ j=1}^n  c_{ij} \left( -\frac{q^{j}_\infty}{2}   + \frac{ \rho^{j}_\infty}{4} + \frac{a}{4}\gamma^{j} \right).
\end{align*} 
Again, a summation w.r.t. $i=1,\dots,n$ directly gives 
the equality of fluxes
\begin{align}\label{eq:sum_q0_halfmoment}
\sum_i q^i_\infty =0\ .
\end{align}
For a uniform node  with equal distribution $c_{ij} = \frac{1}{n-1}, i \neq j$ and $0$ otherwise, 
one obtains using the equality of fluxes
\begin{align*}
\frac{n-2}{2n}q^i_\infty   +\frac{\rho^i_\infty}{4} +\frac{a}{4}\gamma^i  
= \frac{n-2}{2n}q^{j}_\infty   + \frac{ \rho^{j}_\infty}{4} + \frac{a}{4}\gamma^{j} .
\end{align*} 
or the invariance of 
\begin{align}
\label{eq:halfMoment_Couple_Inv1}
\rho_\infty +  \frac{2(n-2)}{n}q_\infty + a\gamma.
\end{align}

Further, integrating the kinetic coupling conditions \eqref{kincoup}
with respect to the positive and negative  half moments, we obtain 
\begin{align*}
\rho^i_+ (0) =  \sum_{j =1}^n  c_{ij} \rho_-^{j} (0)\ .
\end{align*} 
With the half moment approximations \eqref{eq:halfMoment_rqp} this reads
\begin{align*}
  \frac{3 q^i_\infty}{4} +  \frac{\rho^i_\infty}{2}       + \left(\frac{3 a+1}{2}\right) \gamma^i  
  =  \sum_{j=1}^n  c_{ij} \left( - \frac{3 q^{j}_\infty}{4}  +    \frac{1 }{2}  \rho^{j}_\infty     + \left(\frac{3 a-1 }{2}\right) \gamma^{j} \right)\ .
 \end{align*} 
 Summing these conditions for $i=1,\dots,n$ yields 
 \begin{align}\label{eq:sumgamma0}
 \sum_{i=1}^{n}\gamma^i=0\ .
 \end{align}
 Thus, in the case of a uniform node, we  derive another  coupling invariant
  \begin{align}\label{eq:halfMoment_Couple_Inv2}
      \rho_\infty    
      + \frac{3(n-2)}{2n} q_\infty+ \left(3 a+ \frac{n-2}{n}\right) \gamma  \ .
   \end{align} 
  
  Alltogether this yields $2n$  conditions at a node, i.e. \eqref{eq:sum_q0_halfmoment}, \eqref{eq:halfMoment_Couple_Inv1}, \eqref{eq:sumgamma0} and \eqref{eq:halfMoment_Couple_Inv2}. 
  In combination with the conditions at infinity we have $3n$ conditions for $3n$ quantities $\gamma^i, \rho^i_\infty, q^i_\infty$.
   
Note that the invariants \eqref{eq:halfMoment_Couple_Inv2} and  \eqref{eq:halfMoment_Couple_Inv1} can be combined such that $\gamma$ is eliminated, which gives the invariance of
 \begin{align}\label{eq:invariant_halfmoment}
 \rho_\infty + \left(\frac{n-2}{n}\right)\frac{9a+4\left(\frac{n-2}{n}\right)}{4a+2\left(\frac{n-2}{n}\right)}\ q_\infty\ .
  \end{align}

 All the  coupling conditions for the wave equation with uniform nodes derived up to now are given by the conservation of mass \eqref{eq:sum_q_0} and an invariant of the form 
 \begin{align*}
       \rho   + C q\ .
 \end{align*} 
 They differ in the factor $C$, see Table \ref{tableC}.

	\begin{table}
		\begin{center}
	\begin{tabular}{|c|c|}
		\hline
		Coupling conditions& C\\
		\hline
		wave full moment \eqref{eq:invariant_fullmoment}&$\frac{3}{2} \frac{n-2}{n}$\\
		\hline   
		wave Maxwell \eqref{eq:Invariants_Maxwell}& $2 \frac{n-2}{n}$\\
		\hline  
		wave half moment\eqref{eq:invariant_halfmoment}& $\left(\frac{n-2}{n}\right)\frac{\frac{9}{\sqrt{3}}+4\left(\frac{n-2}{n}\right)}{\frac{4}{\sqrt{3}}+2\left(\frac{n-2}{n}\right)}$\\
		\hline        
		\end{tabular}
		\end{center}
		\caption{Coefficients of invariant for different coupling conditions.}
		\label{tableC}
	\end{table}

   With $a =\frac{1 }{\sqrt{3}}$ and $n=3$ the half moment approximation  gives  $C=0.7313$
   compared to $0.6666$ for Maxwell \eqref{eq:Invariants_Maxwell} or $0.5$ for the full moment approximation.
   A numerical comparison of the resulting network solutions is presented in the next section.

\begin{remark}
	For the solution of the wave equation \eqref{euler}, we consider the mathematical entropy
	\begin{align*}
		e = \frac{1}{2}\left(\rho^2+\frac{1}{a^2}q^2\right).
	\end{align*}
	It evolves according to the conservation law
	\begin{align*}
		\partial_t e + \partial_x \left(\rho q\right) = 0 \ .
	\end{align*}
	Along one edge this entropy  is conserved, but the total entropy in the network can change according to the entropy-fluxes at the nodes.
	Note that for all above models with $\rho + Cq= \tilde{C}$ and $C>0$ the total entropy decays, since
	\begin{align*}
		\sum_{i=1}^{n}\rho\, q = 
		\sum_{i=1}^{n}\left(\tilde{C}-Cq\right)\, q = 
		\tilde{C} \sum_{i=1}^{n}q -C\sum_{i=1}^{n}q^2 =  -C\sum_{i=1}^{n}q^2 <0\ .
	\end{align*}	
\end{remark}   
 In the following we discuss also the case of an underlying kinetic model with unbounded velocities.
\subsection{The case of unbounded velocity domain}
We now apply the above procedure to the linear BGK model \eqref{bgk} with an unbounded velocity space, which has \eqref{euler} with arbitrary $a$ as limit equation.

\subsubsection{Half moment equations}
The half moments are defined as
\begin{equation*}
\begin{aligned}
\rho ^- &=\dis\int_{-\infty}^0 f(v)dv\ , \qquad
& \rho^+ &=\dis\int_{0}^\infty f(v)dv\ ,
\\
q ^- &=\dis\int_{-\infty}^0 vf(v)dv \ ,
& q^+ &=\dis\int_{0}^\infty vf(v)dv \ .
\end{aligned}
\end{equation*}
As closure we consider the following approximations of the distribution function  
\begin{align*}
f(v) = (a^+ + vb^+)M(v),\ v\geq 0&&\text{and}&&
f(v) = (a^- + vb^-)M(v),\ v\leq 0
\ ,
\end{align*}
which leads to
\begin{align*}
\begin{aligned}
\rho^- &= \frac{1}{2}a^- -\dis\frac{a}{\sqrt{2\pi}}b^-
\ ,\quad&
\rho^+ &=\frac{1}{2}a^+ + \dis\frac{a}{\sqrt{2\pi}}b^+\ , \\[0.3cm]
q^- &= -\dis\frac{a}{\sqrt{2\pi}}a^- + \dis\frac{a^2}{2}b^- 
\ ,\quad &
q^+ &= \dis\frac{a}{\sqrt{2\pi}}a^+ + \dis\frac{a^2}{2}b^+\ .
\end{aligned}
\end{align*}
Inverting these relations
gives
\begin{align*}
\dis\int_0^\infty v^2f(v)dv &= 
 \frac{1}{\pi - 2}\left((\pi - 4)a^2 \rho^++ a\sqrt{2\pi}q^+\right)
\ ,
\\
\dis\int_{-\infty}^0 v^2f(v)dv &= 
 \frac{1}{\pi - 2}\left((\pi - 4)a^2 \rho^-- a\sqrt{2\pi}q^-\right)
\ .
\end{align*}
This leads to the half-moment system 
\begin{equation*}
\left\{
\begin{array}{lcl}
\partial_t\rho^+ + \partial_xq^+ &=& -\dis\frac{1}{\epsilon}\left(\rho^+ - \left(\dis\frac{\rho }{2} + \frac{ q}{\sqrt{2\pi}a}\right) \right)\\[0.3cm]
\partial_tq^+ + \partial_x\left(\frac{1}{\pi - 2}\left((\pi - 4)a^2 \rho^++ a\sqrt{2\pi}q^+\right)\right) &=& -\dis\frac{1}{\epsilon}\left(q^+ -\left(\dis\frac{a\rho }{\sqrt{2\pi}} + \frac{ q}{2}\right)\right) \\[0.3cm]
\partial_t\rho^- + \partial_xq^- &=& -\dis\frac{1}{\epsilon}\left(\rho^- -
\left(\dis\frac{\rho }{2} -  \frac{ q}{\sqrt{2\pi}a}\right)
\right)\\[0.3cm]
\partial_tq^- + \partial_x\left( \frac{1}{\pi - 2}\left((\pi - 4) a^2\rho^-- a\sqrt{2\pi}q^-\right)\right) &=& -\dis\frac{1}{\epsilon}\left(q^- - \left(- \dis\frac{a\rho }{\sqrt{2\pi}} + \frac{ q}{2}\right)\right) .\\
\end{array}
\right.
\end{equation*}
Introducing the even-odd variables as before 
we can rewrite the system as 
\begin{equation*}
\left\{
\begin{array}{lcl}
\partial_t\rho + \partial_xq &=& 0\\[0.3cm]
\partial_tq + \partial_x\left(\frac{\pi-4}{\pi-2}a^2\rho +  \frac{a\sqrt{2\pi}}{\pi-2}\hat{q}\right) &=&0\\[0.3cm]
\partial_t\hat{\rho} + \partial_x\hat{q} &=&  -\dis\frac{1}{\epsilon}\left(\hat{\rho} -\dis\frac{2}{\sqrt{2\pi} a}q\right)\\[0.3cm]
\partial_t\hat{q} + \partial_x\left(\frac{\pi-4}{\pi-2}a^2\hat{\rho} +\frac{a\sqrt{2\pi}}{\pi-2} q\right) &=& -\dis\frac{1}{\epsilon}\left(\hat{q}-\dis\frac{2a\rho}{\sqrt{2\pi}}\right)\ .
\end{array}
\right.
\end{equation*}


\subsubsection{Half space-half moment problem}

The corresponding half space problem is for 
$x \in \R^+$ given by
\begin{equation}\label{eq:halfspace_halfmoment}
\left\{
\begin{array}{lcl}
\partial_xq &=& 0\\[0.3cm]
\partial_x\left(\frac{\pi-4}{\pi-2}a^2\rho +  \frac{a\sqrt{2\pi}}{\pi-2}\hat{q}\right) &=&0\\[0.3cm]
 \partial_x\hat{q} &=&  -\left(\hat{\rho} -\dis\frac{2}{\sqrt{2\pi} a}q\right)\\[0.3cm]
 \partial_x\left(\frac{\pi-4}{\pi-2}a^2\hat{\rho} +\frac{a\sqrt{2\pi}}{\pi-2} q\right) &=& -\left(\hat{q}-\dis\frac{2a\rho}{\sqrt{2\pi}}\right)\ .
\end{array}
\right.
\end{equation}
As before we need boundary conditions on 
$
\rho^+(0)$, $q^+ (0)$
and a condition at infinity
$$
q_\infty - a \rho_\infty = r_1 = C.
$$
We now construct the explicit solution of this half space problem. 

The first two equations of \eqref{eq:halfspace_halfmoment} state the invariance of 
\begin{align*}
q &= C_1 = const \\
\frac{\pi-4}{\pi-2}a^2\rho +  \frac{a\sqrt{2\pi}}{\pi-2}\hat{q}&= C_2
\ 
\end{align*}
or 
\begin{align*}
\rho = \frac{\pi-2}{a^2(\pi-4)}C_2 - \frac{\sqrt{2\pi}}{a(\pi-4)}\hat{q}.
\end{align*}
From the last equation of \eqref{eq:halfspace_halfmoment} we can deduce that at $x = \infty$ we have
$$
\hat q_\infty = \frac{2a}{\sqrt{2\pi}}\rho_\infty\ ,
$$
which leads with the above invariant to $\rho_\infty = \frac{1}{a^2}C_2 $.
Moreover, from the third line in \eqref{eq:halfspace_halfmoment} we obtain
\begin{align*}
\hat \rho_\infty = \frac{2}{\sqrt{2\pi} a}q_\infty 
\end{align*}
%
%
and thus we can transform the lower part of \eqref{eq:halfspace_halfmoment} to
\begin{equation*}
\left\{
\begin{array}{lcl}
\partial_x\hat{q} &=&  -\left(\hat{\rho} -\dis\frac{2}{\sqrt{2\pi} a}q_\infty\right)\\[0.3cm]
\partial_x\left(\frac{\pi-4}{\pi-2}a^2\hat{\rho} +\frac{a\sqrt{2\pi}}{\pi-2} q_\infty\right) &=& -\left( \hat{q} - \frac{2a\rho}{\sqrt{2\pi}}\left(\frac{\pi-2}{a^2(\pi-4)}C_2 - \frac{\sqrt{2\pi}}{a(\pi-4)}\hat{q}\right)\right)\ .
\end{array}
\right.
\end{equation*}
Rearranging gives
\begin{equation*}
\left\{
\begin{array}{lcrr}
\partial_x\hat{\rho} &=& -\frac{(\pi-2)^2}{(\pi-4)^2a^2} \hat{q} +& \frac{2(\pi-2)^2}{\sqrt{2\pi}a^3(\pi-4)^2}C_2 \\[0.3cm]
\partial_x\hat{q} &=&  -\hat{\rho} +& \frac{2}{\sqrt{2\pi} a}q_\infty
\ .
\end{array}
\right.
\end{equation*}
By defining $ \lambda =  \frac{\pi-2}{a(\pi-4)}$ the solution of $\hat \rho$ and $\hat q $ is
\begin{align*}
\hat \rho= \gamma \lambda \exp(\lambda x)   + \hat \gamma \lambda\exp(-\lambda x) +  \frac{2}{\sqrt{2\pi} a}q_\infty\\
\hat q = -\gamma \exp(\lambda x)   + \hat \gamma \exp( -\lambda x) +\frac{2a}{\sqrt{2\pi}}\rho_\infty .
\end{align*}
Since we are only interested in bounded solutions, we have $\hat \gamma =0$ and thus
\begin{align}
\label{eq:hat_ubounded}
\begin{aligned}
\hat \rho= \gamma \lambda \exp(\lambda x) +  \frac{2}{\sqrt{2\pi} a}q_\infty\\
\hat q = - \gamma \exp(\lambda x)  +\frac{2a}{\sqrt{2\pi}}\rho_\infty .
\end{aligned}
\end{align}
As before, there are three parameters which can be determined with the three conditions
\begin{align*}
\frac{1}{2}\left( \rho (0)+\hat \rho (0) \right) = \rho_+ (0) \\
\frac{1}{2}\left( q(0)+\hat q (0) \right) = q_+ (0)\\
q_\infty -a \rho_\infty = r_1 (0)\ .
\end{align*}
Inserting the expressions \eqref{eq:hat_ubounded} yields 
\begin{align*}
\frac{1}{2}\left(\rho_\infty+\left(\lambda +  \frac{\sqrt{2\pi}}{a(\pi-4)}\right)\gamma  +  \frac{2}{\sqrt{2\pi} a}q_\infty \right) &= \rho_+ (0) \\
\frac{1}{2}\left( q_\infty -\gamma  +\frac{2a}{\sqrt{2\pi}}\rho_\infty \right) &= q_+ (0)
\end{align*}
and for the outgoing quantities $\rho_-(0)$ and $q_-(0)$ we obtain
\begin{align*}
\frac{1}{2}\left(  \rho_\infty+\left(-\lambda + \frac{\sqrt{2\pi}}{a(\pi-4)}\right)\gamma  -  \frac{2}{\sqrt{2\pi} a}q_\infty \right) &= \rho_- (0) \\
\frac{1}{2}\left( q_\infty + \gamma  -\frac{2a}{\sqrt{2\pi}}\rho_\infty \right) &= q_- (0)\ .
\end{align*}

\begin{remark}
The Maxwell approximation  is in this case given by the equation
\begin{align*}
q_+ (0) = q_+(\infty)  = \frac{1}{2} (q_\infty + \hat q_\infty)
= \frac{q_\infty}{2} + \frac{a}{\sqrt{2\pi}}\rho_\infty
\ .
\end{align*} 
The outgoing quantities are computed as
\begin{align*}
q_- (0) = q_-(\infty) =  \frac{ q_\infty}{2} - \frac{a}{\sqrt{2\pi}}\rho_\infty\ ,\\ 
\rho_-(0) 
= \frac{\rho_\infty}{2} -  \frac{1}{\sqrt{2\pi} a}q_\infty
\ .
\end{align*} 
\end{remark}
    
\subsubsection{The extrapolation length}  
In order to estimate the quality of the above approximation, we consider
\begin{align*}
v \partial_x f = - \left(f-\left(\rho + \frac{v}{ a^2}  q \right) M (v)\right) 
\end{align*}
with $q =0$ and $k(v) = v M$, in order to determine $f(\infty,v ) = \lambda_\infty M $.
The Maxwell approximation gives
$ \rho_\infty = \frac{\sqrt{2 \pi}}{2 a}$.
With the half moment method we obtain
\begin{align*}
\frac{1}{2}\left(\rho_\infty+\left(\lambda +  \frac{\sqrt{2\pi}}{a(\pi-4)}\right)\gamma  +  \frac{2}{\sqrt{2\pi} a}q_\infty \right) &= \rho_+ (0)  = \frac{1}{\sqrt{2 \pi}}\\
\frac{1}{2}\left( q_\infty -\gamma  +\frac{2a}{\sqrt{2\pi}}\rho_\infty \right) &= q_+ (0) = \frac{1}{2}\ .
\end{align*}
With $q_\infty =0 $ this results in 
\begin{align*}
\rho_\infty  & = \frac{\pi \sqrt{2 \pi} +  2\pi(1+a)  -2 \sqrt{2 \pi} -8 a}{a(\sqrt{2\pi} \pi + 2 \pi -2 \sqrt{2 \pi} -4)}
\end{align*}
For $a=1$ this gives $\rho_\infty = 1.443$. 
Compared to the very accurate value $1.4371$ obtained with a spectral method \cite{Coron} this yields an error of $ 0.4\%$.
Maxwell gives  $1.2533$, which is an error of $ 12.8\%$.
The variational method \cite{LF81,GK95} gives  $1.4245$ or an error of $0.9\%$.

\subsubsection{Coupling conditions}

A half space half integration of the coupling conditions gives
\begin{align*}
q^i_+ (0) &= -\sum_{ j=1}^n c_{ij}  q_-^{j} (0)
\end{align*}
which yields
\begin{align*}
 q^i_\infty -\gamma^i  +\frac{2a}{\sqrt{2\pi}}\rho^i_\infty &= 
-\sum_{ j=1}^n c_{ij} \left( q^{j}_\infty +\gamma^{j}  -\frac{2a}{\sqrt{2\pi}}\rho^{j}_\infty \right) 
\end{align*}
Summing these equations gives  $\sum_i q^i_\infty =0$.
Further, for a uniform node we observe the invariance of
\begin{align}\label{eq:halfmoment_invariant_q}
 \frac{n-2}{n}q^i_\infty -\gamma^i  +\frac{2a}{\sqrt{2\pi}}\rho^i_\infty &= 
 \frac{n-2}{n}q^{j}_\infty -\gamma^{j}  +\frac{2a}{\sqrt{2\pi}}\rho^{j}_\infty 
 \ .
\end{align}
Moreover, from 
\begin{align*}
\rho^i_+ (0) = \sum_{ j=1}^n c_{ij} \rho_-^{j} (0)
\end{align*} 
one obtains
\begin{align*}
 &\rho^i_\infty+\left(\lambda +  \frac{\sqrt{2\pi}}{a(\pi-4)}\right)\gamma^i  +  \frac{2\ q^i_\infty}{\sqrt{2\pi} a} =
  \sum_{ j=1}^n c_{ij}  \left( \rho^{j}_\infty+\left(-\lambda + \frac{\sqrt{2\pi}}{a(\pi-4)}\right)\gamma^{j}  -  \frac{2\ q^{j}_\infty}{\sqrt{2\pi} a} \right)\,.
  \nonumber
\end{align*} 
Summing these  equations we obtain again $\sum_{i=1}^n \gamma^i=0$.
In the case of a uniform node we can further rearrange, such that
\begin{align}\label{eq:halfmoment_invariant_rho}
\begin{aligned}
\rho^i_\infty+\left(\frac{n-2}{n} + \frac{\sqrt{2\pi}}{\pi-2}\right)&\lambda \gamma^i  +  \frac{2}{\sqrt{2\pi} a}\frac{n-2}{n}q^i_\infty \\
&=
\rho^{j}_\infty+\left(\frac{n-2}{n} +  \frac{\sqrt{2\pi}}{\pi-2}\right)\lambda \gamma^{j}  +  \frac{2}{\sqrt{2\pi} a}\frac{n-2}{n}q^{j}_\infty
\end{aligned}
\end{align} 
holds. 
Together with the conditions at infinity we have again $3n$ conditions for $3 n$ quantities $\gamma^i, \rho^i_\infty, q^i_\infty$.
Combining \eqref{eq:halfmoment_invariant_q} and \eqref{eq:halfmoment_invariant_rho} we can eliminate the $\gamma^i$
and obtain as invariant
\begin{align*} 
\rho   
+
\frac{n-2}{n}\,
\frac{\frac{n-2}{n}(\pi-2)\sqrt{2\pi}+4\pi-8}{\sqrt{2\pi}(\pi-4)+\frac{n-2}{n}(2\pi-4)+2\sqrt{2\pi}}\, \frac{1}{a}\,
q \ .
\end{align*} 
or
\begin{align} 
\label{invunbounded}
\rho   
+
\frac{n-2}{n} \frac{1}{a}\,
\frac{4 + \frac{n-2}{n}\sqrt{2 \pi}}{\sqrt{2 \pi}+2 \frac{n-2}{n}} q
\end{align}

The values of the factor $C$ for the coupling invariants $\rho +C q$ for the present case with unbounded velocities are summarized in Table \ref{tableC2}. 
 
	\begin{table}
		\begin{center}
	\begin{tabular}{|c|c|}
		\hline
		Coupling conditions& C\\
		\hline
		wave full moment \eqref{eq:invariant_fullmomentunbounded}&$\frac{\sqrt{2}}{a \sqrt{\pi}} \frac{n-2}{n}$\\
		\hline   
		wave Maxwell \eqref{maxwellunbounded}& $\frac{\sqrt{\pi}}{a \sqrt{2}} \frac{n-2}{n}$\\
		\hline  
		wave half moment\eqref{invunbounded}& $\frac{n-2}{n} \frac{1}{a}\,
		\frac{4 + \frac{n-2}{n}\sqrt{2 \pi}}{\sqrt{2 \pi}+2 \frac{n-2}{n}}$\\
		\hline        
		\end{tabular}
		\end{center}
		\caption{Coefficients of invariant for different coupling conditions in the unbounded case.}
		\label{tableC2}
	\end{table}			

For $a=1$ and $n=3$ this  factor is approximately $0.5079$.
In contrast, the  Maxwell approximation gives $0.4178$.
For $a^2 =\frac{1}{3}$ the factor is 
$\frac{6 \pi}{6}\sim 0.723$ for Maxwell and approximately $0.8797$ for the above half moment formula.
This has to be compared to the case of bounded velocities discussed in the last section, where the factor has been determined as 
$0.7313$ for the half-moment-  and $0.6666$ for the Maxwell approximation.
Thus, depending on the underlying kinetic model different coupling conditions are obtained for the same macroscopic  equation.

\section{Numerical results}
\label{Numerical results}
In this section we compare the numerical results of the different models on networks.
The solutions of the  kinetic equation \eqref{bgkbounded} are compared with the half-moment approximation \eqref{equ:sc half-moment system} and the macroscopic wave equation \eqref{euler} with different coupling conditions  \eqref{eq:Invariants_Maxwell} \eqref{eq:invariant_fullmoment} \eqref{eq:invariant_halfmoment}
.

The networks are composed of coupled edges, each arc is given by an interval $x \in [0,1]$, which is discretized with $400$ spatial cells if not otherwise stated.
In the kinetic model the velocity domain  $[-1,1]$ is discretized with $400$ cells and we choose $\epsilon = 0.001$ if not otherwise stated.

For the advective part of the equations we use an upwind scheme. 
The source term in the kinetic and half moment equations is approximated with the implicit Euler method.
We note that for the wave equation the upwind scheme yields the exact solution by choosing the CFL number equal to $1$.

At the nodes we consider in addition to the coupling conditions discussed above for comparison a coupling based on the assumption of equal density at the node, see  \cite{VZ09,DZ01,GHKLS11}
\begin{align*}
	\rho^i&=\rho^j \qquad i\neq j,\ i,j=1,2,3\\
	\sum_{i=1}^{3}q_i&=0\ .
\end{align*}

In general, at the outer boundaries  of the network, boundary conditions have to be imposed.
For the kinetic problem the values for the ingoing velocities have to be prescribed as mentioned before. For the half-moment approximation the natural boundary conditions are given by integrating the kinetic conditions and using $\rho_+$ and $q_+$ at left boundaries and 
 $\rho_-$ and $q_-$ at right  boundaries.
%
For the wave equations with full moment,  Maxwell and half-moment conditions
we use the corresponding approximations discussed above to provide boundary values for the macroscopic equations. 

We note that for the boundary values at the  right end of the edges we apply the same procedure as detailed in the previous sections for the left end, reversing  the spatial orientation in the half space problem.
More precisely,  a given kinetic inflow $\ell(v)$ for $v\in[-1,0]$ at the right boundary leads to the following boundary data.
In the  kinetic model we directly impose $f(x=1,v) = \ell(v)$ $v\in[-1,0]$, for the  half moment model the quantities $\rho^-=<\ell(v)>_-$ and $q^-=<v\ell(v)>_-$ are fixed.
In case of the wave equation, the right going Riemann invariant $q+a\rho=C$ is given by the inner data.
The remaining information is given by the following  approximations
	\begin{align*}
		\text{ wave Maxwell}&&& 
		<v\ell(v)>_-= -\frac{\rho}{4}+\frac{1}{2}q\\
		\text{ wave half-moment}&&&
		\left\{
		\begin{array}{rl}
		 <\ell(v)>_-&=\frac{1}{2}\rho - \frac{3}{4}q +\frac{-3a-1}{2}\gamma, \\
		 <v\ell(v)>_-&=-\frac{1}{4}\rho + \frac{1}{2}q +\frac{a}{4}\gamma
		\end{array}\right.
		\\
		\text{ wave full-moment}&&&
		<\ell(v)>_-=\frac{1}{2}\rho - \frac{3}{4}q\ ,
	\end{align*}
	which can be obtained by revisiting section \ref{boundedvel} with $x\in[-\infty,0]$.
	These are used to prescribe the value of the left going characteristic $q-a\rho$.

\subsection{Tripod network}
In a first example, we consider a tripod network with initial conditions
$f^1(x,v) = \frac{1}{2}$,$f^2(x,v) = \frac{1}{3}$ and $f^3(x,v) = 0$.
The corresponding macroscopic states are
$(\rho^1,q^1) = (1,0)$,
$(\rho^2,q^2) = (\frac{2}{3},0)$ and 
$(\rho^3,q^3) = (0,0)$.
We use free boundary conditions at the exterior boundaries.
The computational time is chosen such that the waves generated at the node do not reach the exterior boundaries.

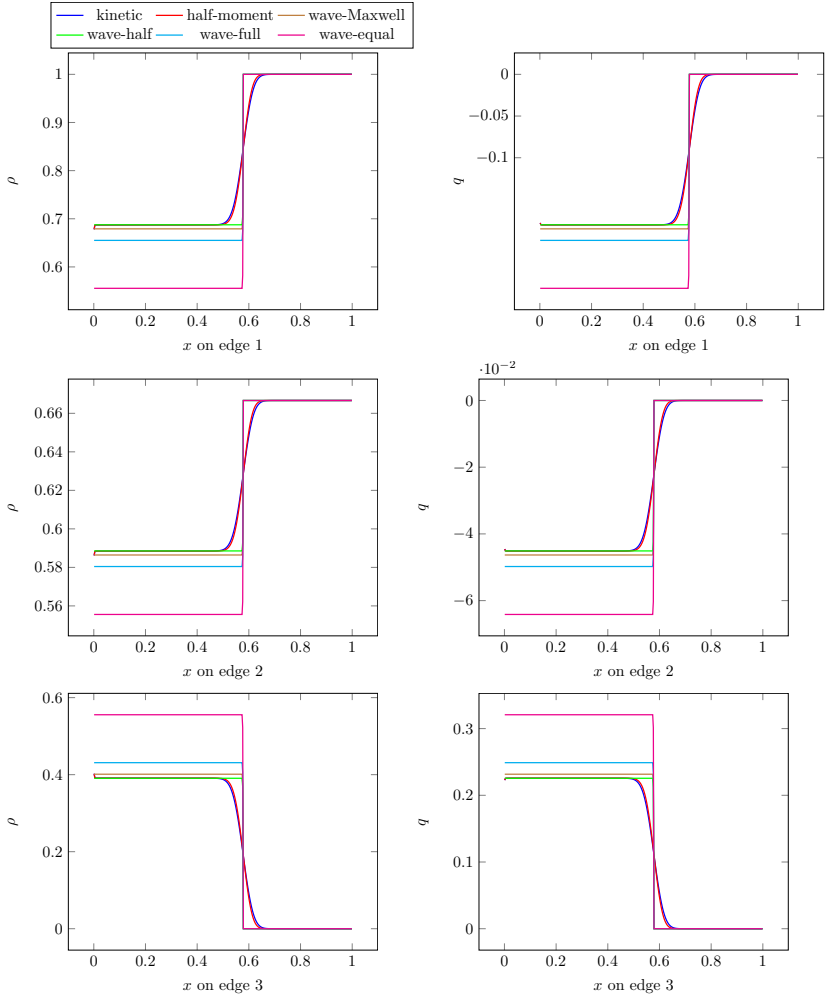
\begin{figure}[ht!]
		\externaltikz{tripod_models_1}{
			\begin{tikzpicture}[scale=0.6]
			\begin{axis}[ylabel = $\rho$,xlabel =  $x$ on edge $1$,
			legend style = {at={(0.5,1)},xshift=0.2cm,yshift=0.1cm,anchor=south},
			legend columns= 3,
			]
			\addplot[color = blue,thick] file{Data/rho_kinetic_1.txt};
			\addlegendentry{kinetic}
			\addplot[color = red,thick] file{Data/rho_half_1.txt};
			\addlegendentry{half-moment}
			\addplot[color = brown,thick] file{Data/rho_wave_Maxwell_1.txt};
			\addlegendentry{wave-Maxwell}
			\addplot[color = green,thick] file{Data/rho_wave_halfmoment_1.txt};
			\addlegendentry{wave-half}
			\addplot[color = cyan,thick] file{Data/rho_wave_full_1.txt};
			\addlegendentry{wave-full}
			\addplot[color = magenta,thick] file{Data/rho_wave_equal_1.txt};
			\addlegendentry{wave-equal}
			\end{axis}
			\end{tikzpicture}
		}
		\externaltikz{tripod_models1_2}{
			\begin{tikzpicture}[scale=0.6]
			\begin{axis}[ylabel = $q$,xlabel = $x$ on edge $1$,
			legend style = {at={(0.5,1)},xshift=0.2cm,yshift=-0.0cm,anchor=south},
			ytick={-0.1,-0.05,0},
			yticklabels={$-0.1$,$-0.05$,$0$},
			]
			\addplot[color = blue,thick] file{Data/q_kinetic_1.txt};
			\addplot[color = red,thick] file{Data/q_half_1.txt};
			\addplot[color = brown,thick] file{Data/q_wave_Maxwell_1.txt};
			\addplot[color = green,thick] file{Data/q_wave_halfmoment_1.txt};
			\addplot[color = cyan,thick] file{Data/q_wave_full_1.txt};
			\addplot[color = magenta,thick] file{Data/q_wave_equal_1.txt};
			\end{axis}
			\end{tikzpicture}
		}
					
		\externaltikz{tripod_models_2}{			
			\begin{tikzpicture}[scale=0.6]
			\begin{axis}[ylabel = $\rho$,xlabel = $x$ on edge $2$,
			legend style = {at={(0.5,1)},xshift=0.2cm,yshift=-0.0cm,anchor=south},
			]
			\addplot[color = blue,thick] file{Data/rho_kinetic_2.txt};
			\addplot[color = red,thick] file{Data/rho_half_2.txt};
			\addplot[color = brown,thick] file{Data/rho_wave_Maxwell_2.txt};
			\addplot[color = green,thick] file{Data/rho_wave_halfmoment_2.txt};
			\addplot[color = cyan,thick] file{Data/rho_wave_full_2.txt};
			\addplot[color = magenta,thick] file{Data/rho_wave_equal_2.txt};
			\end{axis}
			\end{tikzpicture}
		}
		\externaltikz{tripod_models2_2}{
			\begin{tikzpicture}[scale=0.6]
			\begin{axis}[ylabel = $q$,xlabel = $x$ on edge $2$,
			legend style = {at={(0.5,1)},xshift=0.2cm,yshift=-0.0cm,anchor=south},
			]
			\addplot[color = blue,thick] file{Data/q_kinetic_2.txt};
			\addplot[color = red,thick] file{Data/q_half_2.txt};
			\addplot[color = brown,thick] file{Data/q_wave_Maxwell_2.txt};
			\addplot[color = green,thick] file{Data/q_wave_halfmoment_2.txt};
			\addplot[color = cyan,thick] file{Data/q_wave_full_2.txt};
			\addplot[color = magenta,thick] file{Data/q_wave_equal_2.txt};
			\end{axis}
			\end{tikzpicture}
		}	
		
		\externaltikz{tripod_models_3}{
			\begin{tikzpicture}[scale=0.6]
			\begin{axis}[ylabel = $\rho$,xlabel = $x$ on edge $3$,
			legend style = {at={(0.5,1)},xshift=0.2cm,yshift=-0.0cm,anchor=south},
			]
			\addplot[color = blue,thick] file{Data/rho_kinetic_3.txt};
			\addplot[color = red,thick] file{Data/rho_half_3.txt};
			\addplot[color = brown,thick] file{Data/rho_wave_Maxwell_3.txt};
			\addplot[color = green,thick] file{Data/rho_wave_halfmoment_3.txt};
			\addplot[color = cyan,thick] file{Data/rho_wave_full_3.txt};
			\addplot[color = magenta,thick] file{Data/rho_wave_equal_3.txt};
			\end{axis}
			\end{tikzpicture}
		}
		\externaltikz{tripod_models3_2}{
			\begin{tikzpicture}[scale=0.6]
			\begin{axis}[ylabel = $q$,xlabel = $x$ on edge $3$,
			legend style = {at={(0.5,1)},xshift=0.2cm,yshift=-0.0cm,anchor=south},
			]
			\addplot[color = blue,thick] file{Data/q_kinetic_3.txt};
			\addplot[color = red,thick] file{Data/q_half_3.txt};
			\addplot[color = brown,thick] file{Data/q_wave_Maxwell_3.txt};
			\addplot[color = green,thick] file{Data/q_wave_halfmoment_3.txt};
			\addplot[color = cyan,thick] file{Data/q_wave_full_3.txt};
			\addplot[color = magenta,thick] file{Data/q_wave_equal_3.txt};
			\end{axis}
			\end{tikzpicture}
		}
		\caption{Kinetic equation, half-moment equation
		and  wave equation with coupling conditions given by the 
		Maxwell, the half-moment and the  full-moment approach and the assumption of equal density
		at time  $T = 1$.}
		\label{fig1}
		\end{figure}

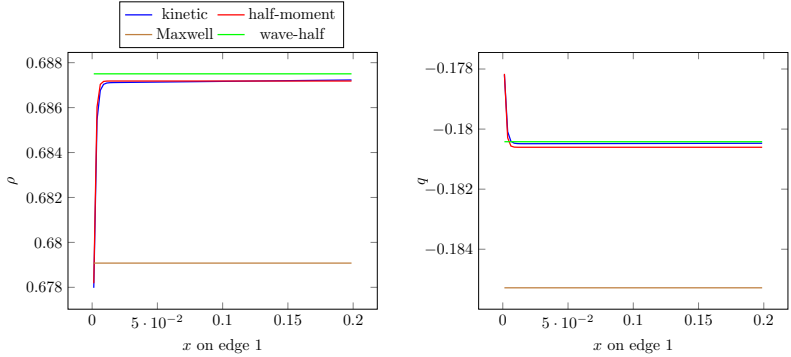
\begin{figure}[ht!]
		\externaltikz{tripod_models_zoom}{			
			\begin{tikzpicture}[scale=0.6]
			\begin{axis}[ylabel = $\rho$,xlabel =  $x$ on edge $1$,
			legend style = {at={(0.5,1)},xshift=0.2cm,yshift=0.1cm,anchor=south},
			legend columns= 2,
			restrict x to domain=0:0.2,restrict y to 
			domain =0.66:0.72,
			yticklabel style={/pgf/number format/fixed,
				/pgf/number format/precision=4}
			]
			\addplot[color = blue,thick] file{Data/rho_kinetic_1.txt};
			\addlegendentry{kinetic}
			\addplot[color = red,thick] file{Data/rho_half_1.txt};
			\addlegendentry{half-moment}
			\addplot[color = brown,thick] file{Data/rho_wave_Maxwell_1.txt};
			\addlegendentry{Maxwell}
			\addplot[color = green,thick] file{Data/rho_wave_halfmoment_1.txt};
			\addlegendentry{wave-half}
			\end{axis}
			\end{tikzpicture}
		}
		\externaltikz{tripod_models_zoom_2}{
			\begin{tikzpicture}[scale=0.6]
			\begin{axis}[ylabel = $q$,xlabel = $x$ on edge $1$,
			legend style = {at={(0.5,1)},xshift=0.2cm,yshift=-0.0cm,anchor=south},	
			restrict x to domain=0:0.2,restrict y to 
			domain =-0.19:-0.17,
			yticklabel style={/pgf/number format/fixed,
			/pgf/number format/precision=4}
			]
			\addplot[color = blue,thick] file{Data/q_kinetic_1.txt};
			\addplot[color = red,thick] file{Data/q_half_1.txt};
			\addplot[color = brown,thick] file{Data/q_wave_Maxwell_1.txt};
			\addplot[color = green,thick] file{Data/q_wave_halfmoment_1.txt};
			\end{axis}
			\end{tikzpicture}
		}
		\caption{Magnification of the situation from Figure \ref{fig1} on edge 1 at the node.}
		\label{fig2}
		\end{figure}

In Figure \ref{fig1} we compare the kinetic, the half-moment
and the wave equation with coupling conditions given by the 
Maxwell, the half-moment and the  full-moment approach and the assumption of equal density
at time  $T = 1$.
	We observe first, that the half moment model gives a very accurate approximation of the kinetic equation.
	Considering the wave equation with different coupling conditions one observes that the interior state is very accurately  approximated by the half moment coupling conditions. 
	Also the Maxwell approximation provides a  good approximation in this case.
	Note that a boundary layer is appearing at the node in the kinetic model, see Figure \ref{fig2} for a magnification of the situation on edge $1$. 
			
	Moreover, we investigate the evolution of the total entropy in the network, i.e.
	$$\sum_{i=1}^{3}\frac{1}{2} \int  \left( (\rho^i)^2+ \frac{1}{a^2}(q^i)^2 \right) dx.$$
	Initially, the total entropy at $t=0$ is equal to $0.722222$.
	In this case we use a very fine grid with $30000$ spatial cells for all models and $400$ cells in  velocity space for the kinetic equation.
	In Table \ref{table1} the value of the total entropy at time $T=0.1$ is shown for the different coupling conditions together with the half moment and the kinetic solution for comparison.
	One observes the very accurate approximation given by the wave equation with half moment coupling conditions.
	Note, that even for the equal density a small amount of  entropy  is lost. 
	This is caused by the numerical diffusion in the very last time step, since we have $\Delta t<\frac{\Delta x}{a}$ to hit the final time.
			\begin{table}[ht!]
				\begin{center}
			\begin{tabular}{|c|c|c|}
				\hline
				Coupling conditions& total entropy & entropy loss\\
				\hline
				wave equal density &$7.2222139e-01 $&$ -8.3353857e-07$\\
				\hline
				wave full moment \eqref{eq:invariant_fullmoment}&$7.1701785e-01$ &$  -5.2043681e-03$\\
				\hline   
				wave Maxwell \eqref{eq:Invariants_Maxwell}&$7.1621400e-01 $&$ -6.0082256e-03 $\\
				\hline  
				wave half moment\eqref{eq:invariant_halfmoment}& $7.1597274e-01 $&$-6.2494790e-03$\\
				\hline        
				half moment \eqref{halfcoupling},$\epsilon = e^{-6}$& $7.1586300e-01 $&$-6.3592218e-03$\\
				\hline   
				kinetic \eqref{kincoup},$\epsilon = e^{-6}$& $7.1574793e-01  $&$-6.4742970e-03$\\
				\hline    
				\end{tabular}
				\end{center}
				\caption{Total entropy and entropy loss at time $T=0.1$ for  different coupling conditions.}
				\label{table1}
			\end{table}			
In Table \ref{tablex} the value of the total entropy at time $T=0.1$ is shown for kinetic and half moment equations for different spatial discretizations. 
			\begin{table}[h!]
				\begin{center}
			\begin{tabular}{|c|c|c|}
				\hline
				 grid &kinetic & half-moment  \\
				\hline        
				 $1000$& $7.1521301e-01 $&$7.1534650e-01  $\\
				\hline   
				 $5000$& $7.1562057e-01 $&$7.1569410e-01  $\\
				\hline   
				 $10000$& $7.1570580e-01  $&$7.1577776e-01 $\\
				\hline   
				$30000$& $7.1574793e-01  $&$7.1586300e-01 $\\
				\hline
				\end{tabular}
				\end{center}
				\caption{Total entropy  for kinetic and half moment for $\epsilon= 10^{-6}$ and different grid size  at time $T=0.1$.}
				\label{tablex}
			\end{table}		

Finally, we investigate the kinetic equation for different values of $\epsilon$.
For  $\epsilon\rightarrow 0$ the kinetic and the half moment model are very well approximated by the solution of the wave equation with half moment coupling conditions, see Figure \ref{fig3} and the corresponding magnification at the node in Figure \ref{fig4}.

	\begin{figure}[ht!]											
		\externaltikz{tripod_epsilon}{
			\begin{tikzpicture}[scale=0.6]
			\begin{axis}[ylabel = $\rho$,xlabel = $x$ on edge $1$,
			legend style = {at={(0.5,1)},xshift=0.2cm,yshift=0.1cm,anchor=south},
			legend columns= 2,
			]
			\addplot[color = blue,thick] file{Data/rho_kinetic_eps01_1.txt};
			\addlegendentry{kinetic $\epsilon = 0.1$}
			\addplot[color = red,thick] file{Data/rho_half_eps01_1.txt};
			\addlegendentry{half-moment $\epsilon = 0.1$}
			\addplot[color = blue,thick] file{Data/rho_kinetic_1.txt};
			\addlegendentry{kinetic $\epsilon = 0.01$}
			\addplot[color = red,thick] file{Data/rho_half_1.txt};
			\addlegendentry{half-moment $\epsilon = 0.01$}
			\addplot[color = blue,thick] file{Data/rho_kinetic_eps0001_1.txt};
			\addlegendentry{kinetic $\epsilon = 0.001$}
			\addplot[color = red,thick] file{Data/rho_half_eps0001_1.txt};
			\addlegendentry{half-moment $\epsilon = 0.001$}
			\addplot[color = green,thick] file{Data/rho_wave_halfmoment_1.txt};
			\addlegendentry{wave half $\epsilon = 0.001$}
			\end{axis}
			\end{tikzpicture}
		}
		\externaltikz{tripod_epsilon2}{
			\begin{tikzpicture}[scale=0.6]
			\begin{axis}[ylabel = $q$,xlabel = $x$ on edge $1$,
			legend style = {at={(0.5,1)},xshift=0.2cm,yshift=-0.0cm,anchor=south},
			]
			\addplot[color = blue,thick] file{Data/q_kinetic_eps01_1.txt};
			\addplot[color = red,thick] file{Data/q_half_eps01_1.txt};
			\addplot[color = blue,thick] file{Data/q_kinetic_1.txt};
			\addplot[color = red,thick] file{Data/q_half_1.txt};
			\addplot[color = blue,thick] file{Data/q_kinetic_eps0001_1.txt};
			\addplot[color = red,thick] file{Data/q_half_eps0001_1.txt};
			\addplot[color = green,thick] file{Data/q_wave_halfmoment_1.txt};
			\end{axis}
			\end{tikzpicture}
		}
			\caption{Solution of kinetic and half moment model for different values of $\epsilon$. For comparison the solution of the wave equation with half moment coupling conditions is shown.}
						\label{fig3}
	\end{figure}
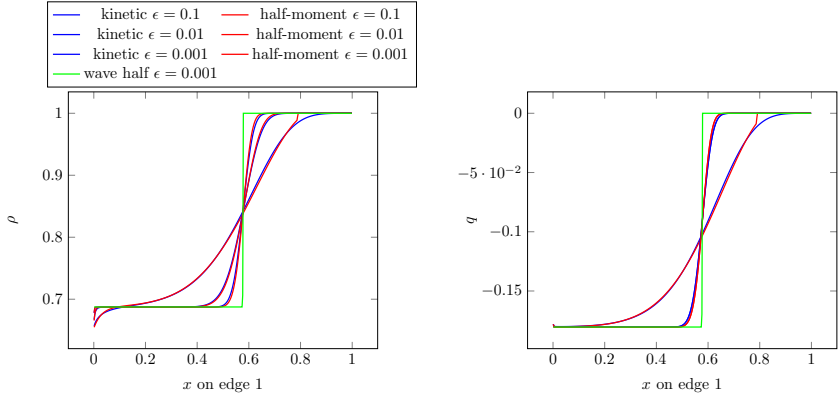

	\begin{figure}[ht!]
		\externaltikz{tripod_epsilon_zoom}{
		\begin{tikzpicture}[scale=0.6]
		\begin{axis}[ylabel = $\rho$,xlabel = $x$ on edge $1$,
		legend style = {at={(0.5,1)},xshift=0.2cm,yshift=-0.0cm,anchor=south},
			legend style = {at={(0.5,1)},xshift=0.2cm,yshift=0.1cm,anchor=south},
			legend columns= 2,
		restrict x to domain=0:0.2,restrict y to domain =0.55:0.7
		]
		\addplot[color = blue,thick] file{Data/rho_kinetic_eps01_1.txt};
		\addlegendentry{kinetic $\epsilon = 0.1$}
		\addplot[color = red,thick] file{Data/rho_half_eps01_1.txt};
		\addlegendentry{half-moment $\epsilon = 0.1$}
		\addplot[color = blue,thick] file{Data/rho_kinetic_1.txt};
		\addlegendentry{kinetic $\epsilon = 0.01$}
		\addplot[color = red,thick] file{Data/rho_half_1.txt};
		\addlegendentry{half-moment $\epsilon = 0.01$}
		\addplot[color = blue,thick] file{Data/rho_kinetic_eps0001_1.txt};
		\addlegendentry{kinetic $\epsilon = 0.001$}
		\addplot[color = red,thick] file{Data/rho_half_eps0001_1.txt};
		\addlegendentry{half-moment $\epsilon = 0.001$}
		\addplot[color = green,thick] file{Data/rho_wave_halfmoment_1.txt};
		\addlegendentry{wave half $\epsilon = 0.001$}
		\end{axis}
		\end{tikzpicture}
		}
		\externaltikz{tripod_epsilon_zoom2}{
		\begin{tikzpicture}[scale=0.6]
		\begin{axis}[ylabel = $q$,xlabel = $x$ on edge $1$,
		legend style = {at={(0.5,1)},xshift=0.2cm,yshift=-0.0cm,anchor=south},
		restrict x to domain=0:0.2,restrict y to domain =-0.2:0.17
		]
		\addplot[color = blue,thick] file{Data/q_kinetic_eps01_1.txt};
		\addplot[color = red,thick] file{Data/q_half_eps01_1.txt};
		\addplot[color = blue,thick] file{Data/q_kinetic_1.txt};
		\addplot[color = red,thick] file{Data/q_half_1.txt};
		\addplot[color = blue,thick] file{Data/q_kinetic_eps0001_1.txt};
		\addplot[color = red,thick] file{Data/q_half_eps0001_1.txt};
		\addplot[color = green,thick] file{Data/q_wave_halfmoment_1.txt};
		\end{axis}
		\end{tikzpicture}
	}
		\caption{Magnification of the situation at the node from Figure \ref{fig3}.}
		\label{fig4}
	\end{figure}
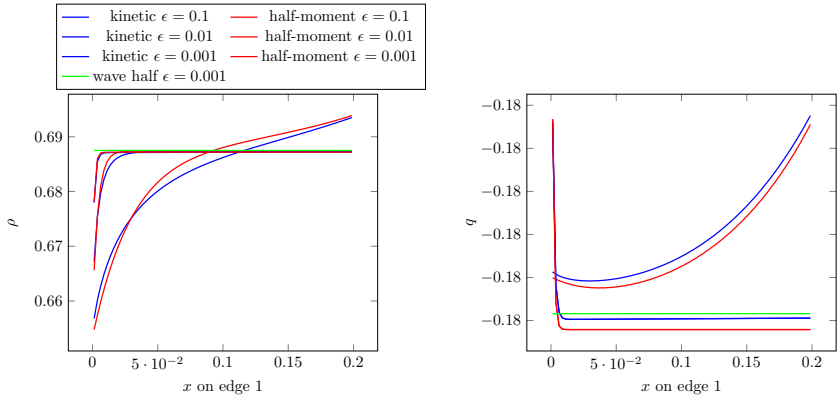
	
\subsection{Diamond network}	
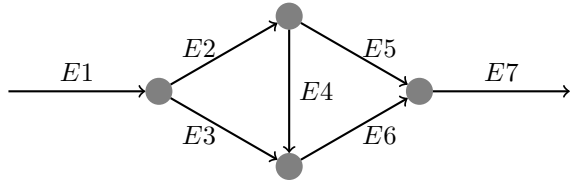
\begin{figure}[ht!]
	\def\len{2}
	\center
	\begin{tikzpicture}[thick,node distance = \len cm]
		\coordinate[] (center);		      
		\node[above of = center,draw,circle,fill,color = gray, yshift=-0.5*\len cm](N2){};
		\node[left of = center,draw,circle,fill,color = gray, xshift=0.134*\len cm](N1){};
		\node[right of = center,draw,circle,fill,color = gray, xshift=-0.134*\len cm](N4){};
		\node[below of = center,draw,circle,fill,color = gray, yshift=0.5*\len cm](N3){};
		\coordinate[left of = N1] (N0);
		\coordinate[right of = N4] (N5);
		\draw[->] (N0)--(N1) node[pos = 0.5,above]{$E1$};
		\draw[->] (N1)--(N2) node[pos = 0.6,left]{$E2\ $};
		\draw[->] (N1)--(N3) node[pos = 0.6,left]{$E3\ $};
		\draw[->] (N2)--(N3) node[pos = 0.5,right]{$E4$};
		\draw[->] (N2)--(N4) node[pos = 0.4,right]{$\ E5$};
		\draw[->] (N3)--(N4) node[pos = 0.4,right]{$\ E6$};
		\draw[->] (N4)--(N5) node[pos = 0.5,above]{$E7$};
	\end{tikzpicture}
	\caption{Diamond network.}
					\label{diamond}
	\end{figure}

As a second example we consider a more complicated network, see Figure \ref{diamond}, as, for example, studied in \cite{EK16} for the wave equation.
		
As initial conditions for the kinetic equation we choose $f^1(x,v) = 1$,$f^2(x,v) = \frac{5}{6}$ and $f^j(x,v) = \frac{1}{2}$ for $j=3,\dots,7$,
which corresponds to macroscopic densities $\rho^1=2$, $\rho^2=\frac{5}{3}$ and $\rho^j=1$ for $j=3,\dots,7$ and fluxes $q^j = 0$ $j=1,\dots,7$.
These data are also prescribed at the two outer boundaries, i.e. $k^1(v)=1$, $v\in [0,1]$ and $\ell^7(v)=\frac{1}{2}$, $v\in [-1,0]$.
Boundary conditions for the wave equation with full moment, Maxwell and half moment conditions are derived as detailed above.
In case of the equal density conditions, we determine the ingoing characteristic using $\rho=1, q=0$ at the $E_1$-boundary and $\rho=\frac{1}{2}, q=0$ at the $E_7$-boundary.

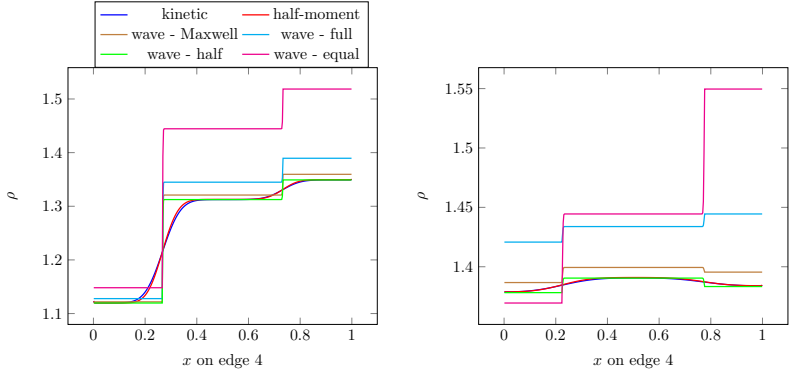
\begin{figure}											
		\externaltikz{edge4}{
			\begin{tikzpicture}[scale=0.6]
			\begin{axis}[ylabel = $\rho$,xlabel = $x$ on edge $4$,
			legend style = {at={(0.5,1)},xshift=0.2cm,yshift=-0.0cm,anchor=south},
			legend columns= 2,
			]
			\addplot[color = blue,thick] file{Data/LargeNetwork/edge4/rho_kinetic_net_T3_4.txt};
			\addlegendentry{kinetic}
			\addplot[color = red,thick] file{Data/LargeNetwork/edge4/rho_half_net_T3_4.txt};
			\addlegendentry{half-moment}
			\addplot[color = brown,thick] file{Data/LargeNetwork/edge4/rho_wave_net_Maxwell_T3_4.txt};
			\addlegendentry{wave - Maxwell}
			\addplot[color = cyan,thick] file{Data/LargeNetwork/edge4/rho_wave_net_full_T3_4.txt};
			\addlegendentry{wave - full}
			\addplot[color = green,thick] file{Data/LargeNetwork/edge4/rho_wave_net_halfmoment_T3_4.txt};
			\addlegendentry{wave - half}
			\addplot[color = magenta,thick] file{Data/LargeNetwork/edge4/rho_wave_net_equal_T3_4.txt};
			\addlegendentry{wave - equal}
			\end{axis}
			\end{tikzpicture}
		}
		\externaltikz{edge4_2}{
			\begin{tikzpicture}[scale=0.6]
			\begin{axis}[ylabel = $\rho$,xlabel = $x$ on edge $4$,
			]
			\addplot[color = blue,thick] file{Data/LargeNetwork/edge4/rho_kinetic_net_T10_4.txt};
			\addplot[color = red,thick] file{Data/LargeNetwork/edge4/rho_half_net_T10_4.txt};
			\addplot[color = brown,thick] file{Data/LargeNetwork/edge4/rho_wave_net_Maxwell_T10_4.txt};
			\addplot[color = cyan,thick] file{Data/LargeNetwork/edge4/rho_wave_net_full_T10_4.txt};
			\addplot[color = green,thick] file{Data/LargeNetwork/edge4/rho_wave_net_halfmoment_T10_4.txt};
			\addplot[color = magenta,thick] file{Data/LargeNetwork/edge4/rho_wave_net_equal_T10_4.txt};
			\end{axis}
			\end{tikzpicture}
		}
		\caption{$\rho$ on edge $4$ at time $t=3$ (left) and time $t=10$ (right).}
		\label{fig5b}
	\end{figure}

In Figure \ref{fig5b} the density $\rho^4$ on edge $4$ is displayed at time $t=3$ and $t=10$.
As before, we observe a good agreement of the half moment coupling with the kinetic and half moment model. 
Also the Maxwell approximation is relatively close to the kinetic results.

\begin{figure}											
		\externaltikz{net_energy}{
			\begin{tikzpicture}[scale=0.6]
			\begin{axis}[ylabel = $\rho$,xlabel = $x$ on edge $4$,
			]
			\addplot[color = blue,thick] file{Data/LargeNetwork/rho_kinetic_net_4.txt};
			\addplot[color = red,thick] file{Data/LargeNetwork/rho_half_net_4.txt};
			\addplot[color = brown,thick] file{Data/LargeNetwork/rho_wave_net_Maxwell_4.txt};
			\addplot[color = cyan,thick] file{Data/LargeNetwork/rho_wave_net_full_4.txt};
			\addplot[color = green,thick] file{Data/LargeNetwork/rho_wave_net_halfmoment_4.txt};
			\addplot[color = magenta,thick] file{Data/LargeNetwork/rho_wave_net_equal_4.txt};
			\end{axis}
			\end{tikzpicture}
		}
		\externaltikz{net_energy_2}{
			\begin{tikzpicture}[scale=0.6]
			\begin{axis}[ylabel = $e$,xlabel = $t$,
			legend style = {at={(0.5,1)},xshift=0.2cm,yshift=0.1cm,anchor=south},
			legend columns= 2,
			]
			\addplot[color = blue,thick] file{Data/LargeNetwork/e_kinetic_net_reduced.txt};
			\addlegendentry{kinetic};
			\addplot[color = red,thick] file{Data/LargeNetwork/e_half_net_reduced.txt};
			\addlegendentry{half-moment}
			\addplot[color = brown,thick] file{Data/LargeNetwork/e_wave_net_reduced_Maxwell.txt};
			\addlegendentry{wave - Maxwell}
			\addplot[color = cyan,thick] file{Data/LargeNetwork/e_wave_net_reduced_full.txt};
			\addlegendentry{wave - full}
			\addplot[color = green,thick] file{Data/LargeNetwork/e_wave_net_reduced_halfmoment.txt};
			\addlegendentry{wave - half}
			\addplot[color = magenta,thick] file{Data/LargeNetwork/e_wave_net_reduced_equal.txt};
			\addlegendentry{wave - equal}
			\end{axis}
			\end{tikzpicture}
		}
		\caption{Left: $\rho$ on edge $4$ at time $t=50$. Right: Evolution of the entropy in the network up to time $t=50$.}
		\label{fig5}
	\end{figure}
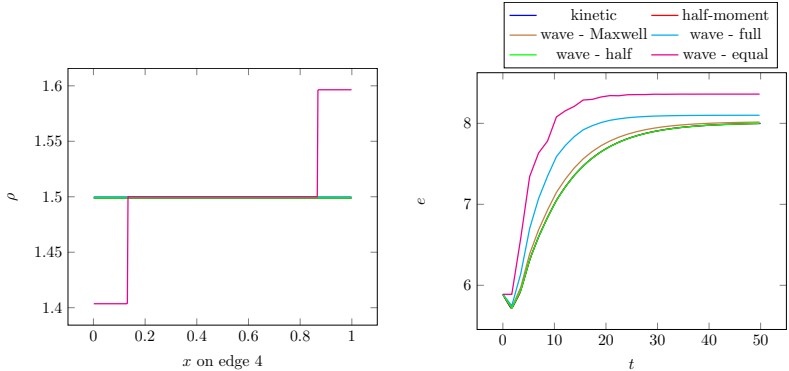

The states of the full moment coupling and the equal density coupling deviate remarkably from the kinetic results.
In Figure \ref{fig5} on the left, $\rho^4$ at time $t=50$ is shown. 
All models have reached a stationary state except the equal density coupling.
Since the entropy is conserved, a set of waves remains trapped in the network, oscillating back and forth.
The stationary states of the models with entropy losses almost coincide.
On the right hand side of Figure \ref{fig5} the evolution of the total entropy over time is plotted. 
The total entropy is increasing due to inflow from the boundaries. 
All models saturate at a certain level, but we again observe a deviation of the equal density and the full moment coupling compared to the  accurate  Maxwell coupling. In this situation the results for the  half-moment coupling, the half moment model and the kinetic model coincide.

\section{Conclusion and Outlook}

In this work we have derived coupling conditions for the wave equation on a network based  
on coupling conditions for an underlying  kinetic BGK type model via a layer analysis  of the situation near the nodes.  
The presentation in this work  includes a new half-moment approximation of  the kinetic  half-space problem.
The general approach can be extended to more  complicated problems like linearized Euler equations or kinetic based coupling conditions for nonlinear problems like Burgers and Lighthill-Whitham type equations which will be considered in a forthcoming publication.

\section*{ Acknowledgment}\ 
\\The first author is supported by the Deutsche \mbox{Forschungsgemeinschaft} (DFG) grant BO 4768/1. The second author by DFG KL 1300/26.
Moreover, funding by the DFG
within the RTG 1932 "Stochastic Models for Innovations in the Engineering
Sciences" is gratefully acknowledged.


\end{document}